\documentclass[12pt,reqno]{amsart}
\usepackage{amsfonts}
\usepackage{natbib}
\usepackage{amssymb}
\usepackage{graphicx}
\usepackage{dsfont}
\usepackage{amsmath}
\usepackage[colorlinks=true,linkcolor=red,citecolor=blue]{hyperref}
\hypersetup{
pdftitle={Dual divergences estimation
for censored survival data},
pdfauthor={\textcopyright\ Cherfi},
}
\usepackage{mathrsfs}
\usepackage{hypernat}
\newcommand{\dint} {\displaystyle\int}

\theoremstyle{definition} \theoremstyle{remark}
\numberwithin{equation}{section}

\renewcommand{\cite}{\citet}
\numberwithin{equation}{section}

\begin{document}

\title{Dual divergences estimation
for censored survival data}
\author{Mohamed Cherfi}
\address{Laboratoire de Statistique Th\'eorique et Appliqu\'ee (LSTA)\\ Equipe d'Accueil 3124\\ Universit\'e Pierre et Marie Curie -- Paris 6\\ Tour 15-25, 2\`eme \'etage\\ 4 place Jussieu\\  75252 Paris cedex 05 e-mail adress :  mohamed.cherfi@gmail.com}
\date{}
\maketitle
\begin{abstract}
This paper is devoted to robust estimation based on dual divergences estimators for parametric models in the framework of right censored data. We give limit laws of the proposed estimators and examine their asymptotic properties through a simulation study.

\end{abstract}

\section{Introduction}
In engineering and biomedical sciences, parametric models are frequently used in analyzing survival data. This analysis is often complicated by the presence of right censoring. Typically right censored data arise in medical studies when patients cannot be followed to the event of interest.

\noindent A common parametric method of estimation is the maximum likelihood approach which is efficient if the specified parametric model is valid. However, in many situations in practice, there is no certainty that the data come
from a specified parametric model and may, in fact, come from some neighborhood of
the model. Likelihood based estimation procedures can lead to poor results when the underlying model is misspecified or contaminated. In such instances, the maximum likelihood is not robust against data or
model inadequacies and the need for robust statistical techniques for estimation and testing has been stressed by
many authors, we may refer to \cite{Huber1981},  \cite{HampelRonchettiRousseeuw1986}, \cite{MaronnaMartinYohai2006} and the references therein.

\noindent In this paper, we consider parametric estimation for right censored
data with and without contamination, and try to balance the dual aims of robustness
and efficiency using minimum divergence estimators.

\cite{Keziou2003} and \cite{BroniatowskiKeziou2009} introduced the class of dual divergences estimators for general parametric models, the procedure being based on the optimization of a new dual form of a divergence and includes the maximum likelihood as a benchmark. \cite{TomaBroniatowski2010} have proved that this class contains robust and efficient estimators and proposed robust test statistics based on divergences estimators.

\noindent A major advantage of the method is that it does not require
additional accessories such as kernel density estimation or other forms of nonparametric
smoothing to produce nonparametric density estimates of the true underlying density function.
The plug-in of the empirical distribution function is sufficient for the purpose of estimating the divergence
in the case of i.i.d. data. For the right-censoring scenario, one can replace the empirical
distribution function with the corresponding estimate of the cumulative distribution function
based on the Kaplan-Meier estimate \cite{KaplanMeier1958}. Thus in this situation one
can also estimate the divergence measure without having to take
recourse to nonparametric smoothing techniques in contrast with existing method, see \cite{Yang1991}, \cite{Ying1992}  that need  a nonparametric estimate of the true density function.  Another feature of the proposed method is it flexibility, that is it leads to a wide class of $M$-estimators indexed by the divergence function and by some instrumental value of the parameter, called here escort parameter. Relevant choices induce efficiency and robustness properties of the proposed estimators.

\noindent The paper is organized as follows. In Section \ref{Section2}, we present the class of dual divergences estimators in the censored case. Asymptotic properties of the proposed estimators are derived in Section \ref{Section3}. We give a brief discussion on the choice of the escort parameter in Section \ref{Section4}. In Section \ref{Section5}, we present Monte Carlo simulation studies to show the performance of the proposed estimators from both robustness and small sample accuracy points of view. Proofs are deferred to the Appendix.

\section{Dual divergences for censored data}\label{Section2}
\noindent The class of dual divergences estimators has been recently introduced by \cite{Keziou2003}, \cite{BroniatowskiKeziou2009}. In the following, we shortly recall their context and definition.

\noindent Recall that the $\phi$-divergence between a bounded
signed measure $Q$ and a probability $P$ on
$\mathscr{D}$, when $Q$ is absolutely continuous with
respect to $P$, is defined by
$$D_\phi(Q,P):=\int_{\mathscr{D}}
\phi\left(\frac{\mathrm{d}Q}{\mathrm{d}P}(x)\right)~\mathrm{d}P(x),$$
where $\phi$ is a convex function from $]-\infty,\infty[$ to
$[0,\infty]$ with $\phi(1)=0$.

Well-known examples of divergences are the Kullback-Leibler, modified Kullback-Leibler, $\chi^2$,
modified $\chi^2$ and {H}ellinger  divergences, they are obtained respectively for
$\phi(x)=x\log x-x+1$, $\phi(x)=-\log x+x-1$,
$\phi(x)=\frac{1}{2}(x-1)^2$,
$\phi(x)=\frac{1}{2}\frac{(x-1)^2}{x}$ and
$\phi(x)=2(\sqrt{x}-1)^2$.
All these divergences belong to the class of the so called
``power divergences'' introduced in \cite{CressieRead1984} (see
also \cite{LieseVajda1987} chapter 2). They are defined through
the class of convex functions
\begin{equation}  \label{powerdivergence}
x\in ]0,+\infty[ \mapsto\phi_{\gamma}(x):=\frac{x^{\gamma}-\gamma
x+\gamma-1}{\gamma(\gamma-1)}
\end{equation}
if $\gamma\in\mathbb{R}\setminus \left\{0,1\right\}$,
$\phi_{0}(x):=-\log x+x-1$ and $\phi_{1}(x):=x\log x-x+1$. (For
all $\gamma\in\mathbb{R}$, we define
$\phi_\gamma(0):=\lim_{x\downarrow 0}\phi_\gamma (x)$). So, the
$KL$-divergence is associated to $\phi_1$, the $KL_m$ to $\phi_0$,
the $\chi^2$ to $\phi_2$, the $\chi^2_m$ to $\phi_{-1}$ and the
{H}ellinger distance to $\phi_{1/2}$. We refer to \cite{LieseVajda1987} for an overview on the
origin of the concept of divergences in statistics.

\noindent Let $X_{1}, \dots, X_{n}$ be an i.i.d. sample
 with p.m. $P_{\theta_{0}}$. Consider the problem of estimating the population parameters of interest $\theta_{0}$, when the underlying identifiable model is given by $\{P_{\theta}: \theta\in \Theta\}$ with $\Theta$ a subset of $\mathbb{R}^{d}$.

\noindent Let $\phi$ be a function of class $\mathcal{C}^2$, strictly convex and satisfies
\begin{equation}
\int \left| \phi ^{\prime }\left( \frac{p_{\theta }(x)}{p_{\alpha
}(x)}\right) \right| ~\mathrm{d}P_{\theta }(x)<\infty .  \label{condition
integrabilite}
\end{equation}
By Lemma 3.2 in \cite{BroniatowskiKeziou2006}, if the function $\phi$ satisfies:
 There exists $0 <\eta < 1$ such that for all $c$ in $\left[1-\eta, 1 + \eta\right]$,
we can find numbers $c_1$, $c_2$, $c_3$ such that
\begin{equation}\label{EqRem1}
\phi(cx)\leq c_1\phi(x) + c_2 \left|x\right| + c_3, \textrm{ for all real } x,
\end{equation}
then the assumption (\ref{condition
integrabilite}) is satisfied whenever $D_{\phi}(P_{\theta}, P_{\alpha})$ is finite. From now on, $\mathcal{U}$ will be the set of $\theta$ and $\alpha$ such that $D_{\phi}(P_{\theta}, P_{\alpha})<\infty$. Note that all the real convex functions $\phi_{\gamma}$ pertaining to the class of power divergences defined in (\ref{powerdivergence}) satisfy the condition (\ref{EqRem1}). Take for example the exponential distribution with density $p_{\theta}(x)=\theta e^{-\theta x}$ for $x\geq 0$ and $\theta>0$, then $\displaystyle{\mathcal{U}:=\left\{\alpha,~\theta>0:~\gamma\theta+(1-\gamma)\alpha>0\right\}}$. 

\noindent Under (\ref{condition integrabilite}), using Fenchel duality technique, the divergence $D_{\phi}(\theta,\theta_0)$ can be represented as resulting from an optimization procedure, this elegant result was proven in \cite{Keziou2003}, \cite{LieseVajda2006} and \cite{BroniatowskiKeziou2009}. \cite{BroniatowskiKeziou2006} called it the dual form of a divergence, due to its connection with convex analysis.

\noindent Under the above conditions, the $\phi$-divergence:
\begin{equation*}
    D_{\phi}(P_{\theta}, P_{\theta_{0}})=\int
\phi\left(\frac{p_{\theta}}{p_{\theta_0}}\right)~\mathrm{d}P_{\theta_0},
\end{equation*}
can be represented as the following form:
\begin{equation}\label{Dualrepresentation}
 D_{\phi}(P_{\theta}, P_{\theta_{0}})=\sup_{\alpha\in
\mathcal{U}}\int h(\theta,\alpha)~\mathrm{d}P_{\theta_{0}},
\end{equation}
where $h(\theta,\alpha):x\mapsto h(\theta,\alpha,x)$ and
\begin{equation}\label{Definition-h}
h(\theta,\alpha,x):=\int \phi ^{\prime }\left( \frac{p_{\theta }}{p_{\alpha
}}\right) ~\mathrm{d}P_{\theta }-\left[ \frac{p_{\theta }(x)}{p_{\alpha
}(x)}\phi ^{\prime }\left( \frac{p_{\theta }(x)}{p_{\alpha }(x)}\right)
-\phi \left( \frac{p_{\theta } (x)}{p_{\alpha }(x)}\right) \right].
\end{equation}
According to \cite{LieseVajda2006}, under the strict convexity and the differentiability of the
function $\phi$, it holds
\begin{equation}\label{Convexity}
\phi(t)\geq\phi(s)+\phi^{\prime}(s)(t-s),
\end{equation}
where the equality holds only for $s=t$. Now, let $\theta$ and $\theta_0$ be fixed
and put $t=p_{\theta}(x)/p_{\theta_0}(x)$ and
$s=p_{\theta}(x)/p_{\alpha}(x)$ in (\ref{Convexity}) and (\ref{Dualrepresentation}) will follow by
integrating with respect to $P_{\theta_0}$.

Since the supremum in (\ref{Dualrepresentation}) is unique and is attained in
$\alpha=\theta_0$, independently upon the value of $\theta$, define the class of estimators of $\theta_{0}$ by
\begin{equation}\label{dualestimator}
\widehat{\alpha}_{\phi}(\theta):=\arg\sup_{\alpha\in \mathcal{U}}\int
h(\theta,\alpha)\mathrm{d}P_{n},\;\;\theta\in \Theta,
\end{equation}
where $h(\theta,\alpha)$ is the function defined in (\ref{Definition-h}). This class is called ``dual
$\phi$-divergence estimators'' (D$\phi$DE's).

\noindent Let us now turn to the estimation using divergences in our setting. In the case of right censored data only
$$Z =\min\left(X, Y\right)\textrm{ and } \delta=1_{\left\{X\leq Y\right\}} $$
are observable. $\delta$ indicates
whether $X$ has been censored or not. The variables
$X_i$ are randomly generated from the true distribution $P_{\theta_{0}}$ which is modeled by the
parametric family $\left\{P_{\theta},~\theta\in\Theta \right\}$. Given a set $\left(Z_i, \delta_i\right),~ i = 1,\ldots,n$ of
independent copies of $\left(Z,\delta\right)$, it is then our goal to draw some inference
on the true but unknown lifetime distribution $P_{\theta_{0}}$.

\noindent Throughout the rest of the paper
we will assume that the variable of interest $X$ and the censoring variable $Y$ are
independent and $G$ denotes the unknown distribution of censoring time $Y$. The distribution $F$ of the observation $Z=\min(X, Y)$, satisfies $1-F=(1-P_{\theta_{0}})(1-G)$.

\noindent \cite{KaplanMeier1958} developed a nonparametric estimator for the survival
function which is is a strongly consistent estimator of the target survival function under appropriate conditions (see \cite{Peterson1977}, \cite{Miller1981})
$$\widehat{P}_{n}(x)=1-\prod_{i=1}^n\left[1-\frac{\delta_{(i)}}{n-i+1}\right]^{\mathds{1}_{\{Z_{(i)\leq x}\}}}$$
where $\left(Z_{(i)}, \delta_{(i)}\right),~ i = 1,\ldots, n$, are the $n$ pairs of observations ordered over
the $Z_{(i)}$ and $\mathds{1}_{A}$ denotes indicator function of $A$.  If  all $\delta_i$'s  are equal to $1$, $\widehat{P}_{n}$ reduces to the ordinary empirical distribution function $P_{n}$.

\noindent Thus, in the right censoring context described above,
we can replace $P_{n}$ in (\ref{dualestimator}) by
$\widehat{P}_{n}(x)$ which provides a consistent
estimator of the true distribution function in this context. Therefore, for the right censoring
situation the ``dual
$\phi$-divergence estimators'' (D$\phi$DE's), is defined by replacing $P_n$ in (\ref{dualestimator}) by $\widehat{P}_{n}$, that is
\begin{equation}\label{dualestimatorcensored}
\widehat{\alpha}_{\phi}(\theta):=\arg\sup_{\alpha\in \mathcal{U}}\int
h(\theta,\alpha)\mathrm{d}\widehat{P}_{n},\;\;\theta\in \Theta.
\end{equation}
\noindent Following \cite{Stute1995}, the {K}aplan-{M}eier integral $\displaystyle{\int
h(\theta,\alpha)\mathrm{d}\widehat{P}_{n}}$  may be written as
\begin{equation*}
    \sum_{i=1}^nW_{in}h(\theta,\alpha,Z_{(i)})
\end{equation*}
where for $1 \leq i \leq n$
\begin{equation*}
  W_{in}=\frac{\delta_{(i)}}{n-i+1}\prod_{j=1}^{i-1}\left[\frac{n-j}{n-j+1}\right]^{\delta_{(j)}}.
\end{equation*}

\noindent The corresponding estimating equation for the unknown parameter is then given by
\begin{equation}\label{estimatingequation}
\int\frac{\partial}{\partial\alpha}h(\theta,\alpha)\mathrm{d}\widehat{P}_{n}=0.
\end{equation}
\noindent Formula (\ref{dualestimatorcensored}) defines a family of $M$-estimators for censored
data indexed by the function $\phi$ specifying the divergence and by some instrumental value of the parameter $\theta$, called here escort parameter, see also \cite{BroniatowskiVajda2009}. The choices of $\phi$ and $\theta$ represent a major feature of the estimation procedure, since they induce efficiency and robustness properties.

\noindent An $M$-estimator of $\psi$-type is the solution of the vector equation:
\begin{equation}\label{general-m-estimator}
    \int\psi(x;\alpha)\mathrm{d}\widehat{P}_n=0,
\end{equation}
where the elements of $\psi(x;\alpha)$ represent the partial derivatives of $h(\theta,\alpha,x)$ with respect to the components of $\alpha$.

\noindent The first extension of $M$-estimators to censored data was  noted in \cite{Reid1981}, she derived the influence function and then the asymptotic normality.  \cite{Oakes1986} considered $M$-estimators (\ref{general-m-estimator}) with $\psi(x;\theta) = -\log f(x;\theta)$ and called them approximate MLEs (hereafter AMLE). \cite{Wang1995} studied the strong consistency of $M$-estimators using the law of large numbers of the {K}aplan-{M}eier integral developed by \cite{StuteWang1993} and \cite{Stute1995}. \cite{Wang1999} extended asymptotic results for $M$-estimators to the censored case.

\noindent The {H}ellinger distance have been used by \cite{Yang1991} and \cite{Ying1992}. Estimation under misspecification have been considered by \cite{SuzukawaImaiSato2001}. \cite{BasuBasuJones2006} developed a robust estimation, adapting the robust density power divergence methodology of \cite{BHHJ1998}.

\section{Asymptotic properties}\label{Section3}
\noindent In this section, we establish the consistency and asymptotic normality of the class of dual divergences estimators in the right censored situation.

\noindent For a distribution $P$, let $\tau_{P}=\sup\left\{x:~ P(x) < 1\right\}$ denote the upper bound of
the support of $P$.

\noindent Assume that $\theta_{0}$ is an interior point of $\Theta$, the convex function $\phi$ has continuous derivatives up to 4th order and the density $p_{\alpha}(x)$ has continuous partial derivatives
up to 3th order (for all $x$ $\lambda-a.e$). Hereafter, $\dot{p}_\alpha$ will denotes the derivative with respect to $\alpha$ of $p_{\alpha}$, $\|\cdot\|$ the Euclidean norm, and, for a real valued function $g$, its total
variation or variation norm is defined as
$$\left\|g\right\|_{\mathbf{v}}=\sup\sum_{j=1}^{N+1}\left|g(x_j)-g(x_{j-1})\right|,$$
where the supremum is taken over all $N$ and over all choices of $\left\{x_j\right\}$ such that
$$-\infty=x_0< x_1 <\ldots<x_N<x_{N+1}=+\infty .$$

\noindent Let $S$ be the $d\times d$ matrix with entries
\begin{equation*}\label{S_generique}
   S_{ij}=-P_{\theta_{0}}\frac{\partial^2}{\partial\alpha_i\partial\alpha_j}h(\theta,\theta_0).
\end{equation*}
\noindent We precise some notations for the asymptotic results in this section. The following quantities have been introduced in \cite{Stute1995b} and \cite{Wang1999}.

\noindent Denote $m(y)=p(\delta=1|Y =y)$, decompose $F$ into two subdistributions $F_0,~F_1$, such that $F=F_0+F_1$, where
\begin{eqnarray*}
\nonumber
  F_0(y) &=&P(Y\leq y,\delta=0)=\int_{-\infty}^y(1-m(t))\mathrm{d}F(t)=\int_{-\infty}^y(1-P_{\theta_0}(t))\mathrm{d}G(t),\\
  F_1(y) &=&P(Y\leq y,\delta=1)=\int_{-\infty}^ym(t)\mathrm{d}F(t)=\int_{-\infty}^y(1-G(t-))~\mathrm{d}P_{\theta_0},
\end{eqnarray*}
and their empirical counterparts
\begin{equation*}
F_{jn}(y)=\frac{1}{n}\sum_{i=1}^n1_{\{Z_i\leq y,\delta_i=j\}},\quad j=0,~1.
\end{equation*}

\noindent Define
\begin{equation}\label{xi0}
   \xi_{0}(x)=\exp\left\{\int\frac{1_{\{y<x\}}\mathrm{d}F_{0}(y)}{1-F(y)}\right\},
\end{equation}
and, for $i=1,\ldots,d$,
\begin{equation}\label{xi1}
   \xi_{1i}(x)=\left[1-F(x)\right]^{-1}\int1_{\{x<y\}}\frac{\partial}{\partial\alpha_i}h(\theta,\alpha,y)\xi_{0}(y)\mathrm{d}F_{1}(y),
\end{equation}
\begin{equation}\label{xi2}
   \xi_{2i}(x)=\int\frac{\partial}{\partial\alpha_i}h(\theta,\alpha,z)\xi_{0}(z)\boldsymbol{C}(x\wedge z)\mathrm{d}F_{1}(z),
\end{equation}
where
\begin{equation}\label{varphi}
\boldsymbol{C}(x)=\int\frac{1_{\{y<x\}}\mathrm{d}F_{0}(y)}{\left[1-F(y)\right]^2}=\int\frac{1_{\{y<x\}}\mathrm{d}G(y)}{\left[1-P_{\theta_0}(y)\right]\left[1-G(y)\right]^2}.
\end{equation}

\noindent Let $U(\alpha)=\left(U_1,\ldots,U_d\right)^{\top}$ denote the random variable defined as:
\begin{equation}\label{Ui}
    U_i(\alpha)=\frac{\partial}{\partial\alpha_i}h(\theta,\alpha,Y)\xi_{0}(Y)\delta+\xi_{1i}(Y)(1-\delta)-\xi_{2i}(Y),~i=1\ldots,d.
\end{equation}
When $\alpha=\theta_0$,
$$U_i(\theta_0)=\frac{\partial}{\partial\alpha_i}h(\theta,\theta_0,Y)\xi_{0}(Y)\delta+\xi_{1i}(Y)(1-\delta)-\xi_{2i}(Y),~i=1\ldots,d.$$
Denote $V$ the $d\times d$ matrix
\begin{equation}\label{Mcensored}
    V=E\left(U(\theta_0)U(\theta_0)^{\top}\right).
\end{equation}
\subsection{Consistency}
\noindent In Theorem \ref{Theo1} below, we prove that $\widehat{\alpha}_\phi(\theta)$ exist and are consistent. We will consider the following conditions.
\begin{enumerate}
\item[(R.0)] $\tau_{P_{\theta_{0}}}\leq \tau_{G}$, where equality may hold except when $G$ is
continuous at $\tau_{P_{\theta_{0}}}$, and,  the probability mass of $P_{\theta_{0}}$ at $\tau_{P_{\theta_{0}}}$: $P_{\theta_{0}}\left(\tau_{P_{\theta_{0}}}\right)-P_{\theta_{0}}\left(\tau_{P_{\theta_{0}}}-\right)> 0$;
\item [(R.1)] There exists a neighborhood $N(\theta_0)$
of $\theta_0$ such that the first and second order partial
derivatives (w.r.t $\alpha$) of $\phi'\left(p_\theta (x)/p_\alpha (x)\right)p_\theta (x)$
are dominated on $N(\theta_0)$ by some integrable
functions. The third order partial derivatives (w.r.t $\alpha$) of
$h(\theta,\alpha,x)$ are dominated on $N(\theta_0)$ by some
$P_{\theta_0}$-integrable functions and the matrices
$S$ and $V$ are non singular;
\item[(R.2)] $\displaystyle{\left\|\frac{\partial}{\partial\alpha}\,
h(\theta ,\theta_0)\right\|_{\mathbf{v}}<\infty}$.
\end{enumerate}
These conditions are mild and can be satisfied in most of circumstances. The condition {\rm{(R.0)}} ensures that $X$ is observable on the hole of the support of $P_{\theta_0}$. Note that if  $\tau_{P_{\theta_{0}}}> \tau_{G}$ holds, the $X_i$ in $[\tau_{G},\infty)$ is certainly censored. In a large number of practical situations, $\tau_{P_{\theta_{0}}}= \tau_{G}=\infty$, hence the condition {\rm{(R.0)}} is satisfied.

\noindent Condition {\rm{(R.1)}} is about usual regularity properties of the underlying model, it guarantees that we can interchange
integration and differentiation and the existence of the variance-covariance matrices, it is similar to regularity conditions used in \cite{Keziou2003} and \cite{BroniatowskiKeziou2009} in the uncensored case.

\noindent Condition {\rm{(R.2)}} is needed to apply the L.I.L in the proof of Theorem \ref{Theo1}. The requirement that $\displaystyle{\psi(x;\alpha):=\frac{\partial}{\partial\alpha}\,
h(\theta ,\alpha)}$ be of bounded variation is standard in $M$-estimation, see for instance \cite{Welsh1989}. Keep in mind the assumed regularity conditions on the criterion function, that is,
$h(\theta,\alpha)$ in the  present framework, to see that it holds for most regular models. 

It is also noted that conditions {\rm{(R.1)}} and {\rm{(R.2)}} are independent of $G$.

\newtheorem{Theo1}{Theorem}
\begin{Theo1}\label{Theo1}{\rm
Let $B(\theta_{0},n^{-1/3}):=\left\{\theta \in
\Theta, \|\theta -\theta_{0} \| \leq
n^{-1/3} \right\}$. Assume that conditions (R.0-2) hold, then as $n$ tends to infinity, with
probability one, the function $\displaystyle{\alpha\mapsto\int
h(\theta ,\alpha )~\mathrm{d}\widehat{P}_{n}}$ attains
its local maximum at some point $\widehat{\alpha}_\phi(\theta)$ in the interior of  $B(\theta_{0},n^{-1/3})$, which implies that the estimate
 $\widehat{\alpha}_\phi(\theta)$ is consistent and satisfies
 \begin{equation*}
\int\frac{\partial}{\partial\alpha}h(\theta,\widehat{\alpha}_\phi(\theta))~\mathrm{d}\widehat{P}_{n}= 0.
 \end{equation*}
 }
 \end{Theo1}
 The proof of Theorem \ref{Theo1} is postponed to the Appendix.

In practice, to obtain the estimate $\widehat{\alpha}_\phi(\theta)$, we use gradient descent algorithms in the optimization in (\ref{estimatingequation}). These algorithms depend on some initial parameter value of $\alpha$. Hence, it is desirable to prove that in a neighborhood of $\theta_0$ there exists a maximum of $\displaystyle{\int h(\theta ,\alpha )~\mathrm{d}\widehat{P}_{n}}$ which does indeed converge to $\theta_0$. Note that the initial parameter value may provide a local maximum (not necessarily global) of $\displaystyle{\int
h(\theta ,\alpha )~\mathrm{d}\widehat{P}_{n}}$.  The local and global estimates coincide if the function $\displaystyle{\alpha\in\Theta\mapsto\int
h(\theta ,\alpha )~\mathrm{d}\widehat{P}_{n}}$ is strictly concave and $\Theta$ is convex, see for instance \cite[Remark 3.5]{BroniatowskiKeziou2009}. 

The aim of Theorem \ref{Theo1} is not to establish the optimal rate
of the estimate but merely the existence and the consistency (a.s.) of the estimate.
We have considered $n^{-1/3}$ because it works well, indeed, in Taylor expansion (\ref{Taylor}),
in the proof, the third term of the right hand side is $O(1)$ only for this rate, which is the major
key of the demonstration, for similar arguments in the estimation of copula models see \cite{BouzebdaKeziou2010}.

\subsection{Asymptotic normality}

\noindent In Theorem \ref{Theo2} below, we  give the limit law of the estimates
$\widehat{\alpha}_\phi(\theta)$ under the following conditions.  From now on, $\overset{d}{\longrightarrow}$ denotes the convergence in distribution. \begin{enumerate}
 \item [(R.3)] For all $1\leq i\leq d$, $\displaystyle{E\left[\left(\frac{\partial}{\partial\alpha_i}h(\theta ,\alpha,Y )\xi_{0}(Y)\delta\right)^2\right]<\infty;}$
 \item [(R.4)] For all $1\leq i\leq d$, $\displaystyle{\int\left|\frac{\partial}{\partial\alpha_i}h(\theta ,\alpha,x )\right|\boldsymbol{C}^{1/2}(x)\mathrm{d}P_{\theta_0}<\infty.}$
\end{enumerate}
Conditions {\rm(R.3-4)} are essential for the asymptotic results of $M$-estimators in the censored case, see for instance \cite{Wang1999} and \cite{BasuBasuJones2006} in the case of density power divergence method.
\newtheorem{Theo2}[Theo1]{Theorem}
\begin{Theo2}\label{Theo2}{\rm
Assume that assumptions (R.0-4) hold. Then, as $n\rightarrow\infty$ $$\sqrt{n}\left(\widehat{\alpha}_\phi(\theta)-\theta_{0}\right)\overset{d}{\longrightarrow}\mathcal{N}\left(0,S^{-1}VS^{-1}\right)$$
}
\end{Theo2}
The proof of Theorem \ref{Theo2} is postponed to the Appendix.

%

\section{Adaptive choice of the escort parameter}\label{Section4}
\noindent Analogously as in the uncensored case, the very peculiar choice of the escort parameter defined through $\theta=\theta_0$ has same limit properties as the AMLE. The D$\phi$DE $\widehat{\alpha}_\phi\left(\theta_{0}\right)$, in this case, has variance which indeed coincides with the AMLE for censored data. If $\theta$ is a real parameter, the asymptotic distribution of $\sqrt{n}\left(\widehat{\alpha}_\phi(\theta)-\theta_{0}\right)$ is normal with mean zero and variance
\begin{equation}
\dfrac{\dint\dfrac{\dot{p}^2_{\theta_0}(x)}{p_{\theta_0}(x)\overline{G}(x)}~\mathrm{d}x-\dint\dfrac{\dot{P}^2_{\theta_0}(x)}{\overline{P}_{\theta_0}(x)\overline{G}^2(x)}~\mathrm{d}x}{I^2_{\theta_0}},
\end{equation}
where $\dot{p}_{\theta}$ is the derivative with respect to
$\theta$ of $p_{\theta}$ and $I_{\theta_0}$ is
the Fisher information matrix
\begin{equation*}
I_{\theta_0}:=\int
\frac{\dot{p}_{\theta_0}\dot{p}_{\theta_0}^{\top}}{
p_{\theta_0}} ~ \mathrm{d}\lambda.
\end{equation*}
Observe that if there is no censorship, that is $G\equiv 0$, the variance of $\widehat{\alpha}_\phi\left(\theta_{0}\right)$ is $\dfrac{1}{I_{\theta_0}}$. 

This result is of some relevance, since it leaves open the choice of the divergence, while keeping good asymptotic properties.

In practice, the consequence is that the escort parameter should be chosen as a the AML estimator of $\theta_0$, say $\widehat{\theta}_{n}$, which under the model is a consistent estimate of $\theta_{0}$. In turn we may expect that the resulting estimator $\widehat{\alpha}_\phi\left(\widehat{\theta}_{n}\right)$ inherits both good asymptotic properties under the model, and, under contamination through a tuning of the divergence index $\gamma$.

Consider the power divergences family \cite{CressieRead1984}, the estimating equation (\ref{estimatingequation}) reduces to
\begin{equation}\label{estimatingequation-pd}
-\int
\left(\frac{p_{\theta}(x)}{p_{\alpha}(x)}\right)^{\gamma-1}\frac{\dot{p}_{\alpha}(x)}{p_{\alpha}(x)}p_{\theta}(x)~\mathrm{d}x+\frac{1}{n}\sum_{i=1}^{n}W_{in}\left(\frac{p_{\theta}(Z_{(i)})}{p_{\alpha}(Z_{(i)})}\right)^{\gamma}\frac{\dot{p}_{\alpha}(Z_{(i)})}{p_{\alpha}(Z_{(i)})}=0,
\end{equation}
where $W_{in}$ are the {K}aplan-{M}eier weights. The estimate $\widehat{\alpha}_\phi(\theta)$ is the solution in $\alpha$ of (\ref{estimatingequation-pd}).

\noindent An improvement of the present estimate results in the plugging of a preliminary consistent estimate of $\theta_0$, say $\widehat{\theta}_{n}$, as an adaptive escort parameter $\theta$ choice.
\begin{figure}[h]
\centerline{\includegraphics[width=8cm]{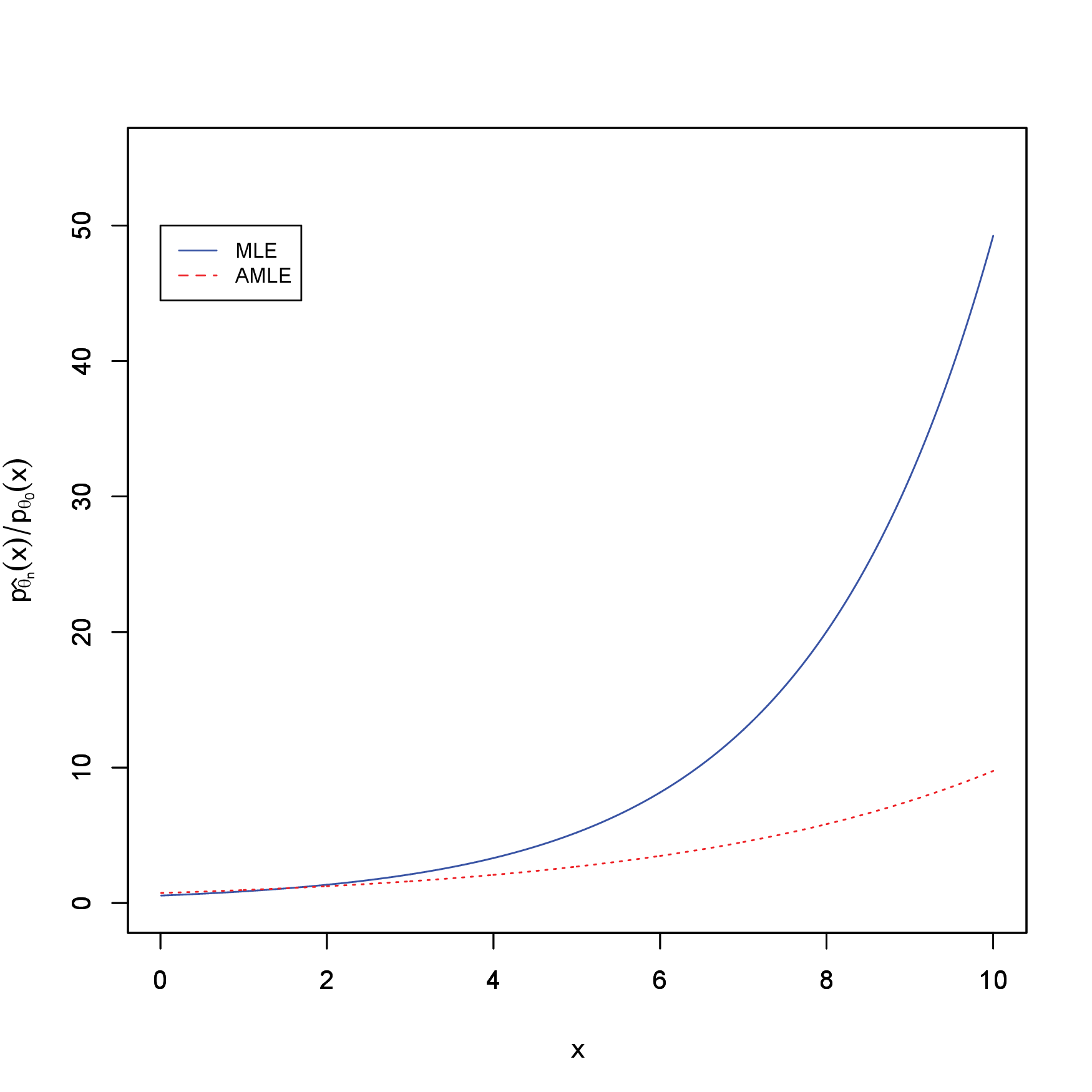}}
\caption{Behaviour of the ratio $\displaystyle{\frac{p_{\widehat{\theta}_{n}}(x)}{p_{\theta_0}(x)}}$ under conatmination, for a randomly generated exponential sample $\exp(1)$ of size $100$  with $\exp(1/9)$ as censoring distribution and $20\%$ of  contamination  by $\exp(0.1)$.} \label{ratio}
\end{figure}

\noindent Let $x$ be some outlier, the role of the outlier $x$ in (\ref{estimatingequation-pd}) appears in the term
\begin{equation}\label{leadingterm}
\left(\frac{p_{\widehat{\theta}_{n}}(x)}{p_{\alpha}(x)}\right)^{\gamma}\frac{\dot{p}_{\alpha}(x)}{p_{\alpha}(x)}.
\end{equation}
The estimate $\widehat{\alpha}_\phi(\theta)$ is robust if this term is stable. That is, if it is small when $\alpha$ is near $\theta_0$. If the escort parameter $\widehat{\theta}_{n}$ is not a robust estimator, the ratio $\displaystyle{\frac{p_{\widehat{\theta}_{n}}(x)}{p_{\theta_0}(x)}}$ can be very large, see Figure \ref{ratio}. This is due to the fact that the outlier $x$ will be more likely under $P_{\widehat{\theta}_{n}}$, that is $\widehat{\theta}_{n}$ will lead to an over evaluation of $p_{\widehat{\theta}_{n}}(x)$ with respect to the expected value under $\theta_0$, say $p_{\theta_0}(x)$.  To guard against such situations, compensate through the choice of $\gamma$, this requires further investigation.

One proposal for the choice of the divergence, is to look for values of the tuning parameter $\gamma$ to obtain a bounded influence function in the spirit of \cite{TomaBroniatowski2010}, we leave this issue open for future research.

We now prove that the subsequent estimator $\widehat{\alpha}_\phi\left(\widehat{\theta}_{n}\right)$ enjoys a limit normal law under the model, see Theorem \ref{Theo3} below.

\noindent Recall that, when $\theta=\theta_0$, $S=-\phi^{\prime\prime}(1)I_{\theta_0}$.
Also, when $\alpha=\theta=\theta_0$, we have
\begin{equation*}
    U=\phi^{\prime\prime}(1)\frac{\dot{p}_{\theta_0}}{p_{\theta_0}}\xi_{0}(Y)\delta+\xi_{1}(Y)(1-\delta)-\xi_{2}(Y),
\end{equation*}
 and the matrix $V$ defined in (\ref{Mcensored}) is
\begin{equation}\label{Vnew}
    V=E\left(UU^{\top}\right),
\end{equation}

\begin{enumerate}
\item [(R.5)] For all $1\leq i,j\leq d$, any one of the following conditions holds:
    \begin{itemize}
  \item[(i)] $\displaystyle{\theta\mapsto\frac{\partial^2}{\partial\alpha_i\partial\alpha_j}h(\theta ,\theta_0,x )}$ is continuous at $\theta_{0}$ uniformly in $x$;
  \item [(ii)] $\displaystyle{\int\sup_{\{\theta:|\theta-\theta_{0}|\leq\rho\}}\left|\frac{\partial^2}{\partial\alpha_i\partial\alpha_j}h(\theta ,\theta_{0} )-\frac{\partial^2}{\partial\alpha_i\partial\alpha_j}h(\theta_{0} ,\theta_{0} )\right|~\mathrm{d}P_{\theta_{0}}=\epsilon_{\rho}\rightarrow 0}$,\\ as $\rho \rightarrow 0$.
  \item[(iii)] $\displaystyle{x\mapsto\frac{\partial^2}{\partial\alpha_i\partial\alpha_j}h(\theta ,\theta_{0},x )}$ is continuous in $x$ for $\theta$ in a neighborhood of $\theta_{0}$ and $$\lim_{\theta\rightarrow\theta_0}\left\|\frac{\partial^2}{\partial\alpha_i\partial\alpha_j}h(\theta ,\theta_{0},\cdot )-\frac{\partial^2}{\partial\alpha_i\partial\alpha_j}h(\theta_{0} ,\theta_0,\cdot )\right\|_{\mathbf{v}}=0;$$
  \item[(iv)]$\displaystyle{\theta\mapsto\int\frac{\partial^2}{\partial\alpha_i\partial\alpha_j}h(\theta ,\theta_{0})\mathrm{d}P_{\theta_0}}$ is continuous at $\theta=\theta_0$, and \\ $\displaystyle{x\mapsto\frac{\partial^2}{\partial\alpha_i\partial\alpha_j}h(\theta ,\theta_{0},x)}$ is continuous in $x$ for $\theta$ in a neighborhood of $\theta_{0}$ and $\displaystyle{\lim_{\theta\rightarrow\theta_0}\left\|\frac{\partial^2}{\partial\alpha_i\partial\alpha_j}h(\theta ,\theta_{0},\cdot )-\frac{\partial^2}{\partial\alpha_i\partial\alpha_j}h(\theta_{0} ,\theta_0,\cdot )\right\|_{\mathbf{v}}<\infty;}$
  \item [(v)]$\displaystyle{\theta\mapsto\int\frac{\partial^2}{\partial\alpha_i\partial\alpha_j}h(\theta ,\theta_{0})\mathrm{d}P_{\theta_0}}$ is continuous at $\theta=\theta_0$, and $$\int\frac{\partial^2}{\partial\alpha_i\partial\alpha_j}h(\theta ,\theta_{0})\mathrm{d}\widehat{P}_{n}\overset{P}{\longrightarrow}\int\frac{\partial^2}{\partial\alpha_i\partial\alpha_j}h(\theta ,\theta_{0})\mathrm{d}P_{\theta_0}<\infty,$$ uniformly for $\theta$ in a neighborhood of $\theta_{0}$.
\end{itemize}
\end{enumerate}

Condition (R.5) is related to Lemma 1 in \cite{Wang1999} and ensures the convergence
 $$\int\frac{\partial^2}{\partial\alpha_i\partial\alpha_j}h(\widehat{\theta}_{n} ,\theta_0)\mathrm{d}\widehat{P}_{n}\overset{P}{\longrightarrow}\int\frac{\partial^2}{\partial\alpha_i\partial\alpha_j}h(\theta_0 ,\theta_0)~\mathrm{d}P_{\theta_0},~1\leq i,j\leq d,$$
 provided that $\displaystyle{\int\left|\frac{\partial^2}{\partial\alpha_i\partial\alpha_j}h(\theta_0 ,\theta_0)\right|\mathrm{d}P_{\theta_0}<\infty,~1\leq i,j\leq d}$, $\widehat{\theta}_{n}\overset{P}{\longrightarrow}\theta_{0}$  and condition (R.0) holds.

\newtheorem{Theo3}[Theo1]{Theorem}
\begin{Theo3}\label{Theo3}
{\rm
Assume that assumptions (R.0-5) hold. Then, as $n\rightarrow\infty$ $$\sqrt{n}\left(\widehat{\alpha}_\phi\left(\widehat{\theta}_{n}\right)-\theta_{0}\right)\overset{d}{\longrightarrow}\mathcal{N}\left(0,\phi^{\prime\prime-2}(1)I_{\theta_0}^{-1}VI_{\theta_0}^{-1}\right),$$
where $V$ is defined in (\ref{Vnew}).
}
\end{Theo3}
The proof of Theorem \ref{Theo3} is postponed to the Appendix.

\section{Simulation}\label{Section5}
In this section, we present results of a simulation study which was conducted
to explore the properties of newly proposed dual $\phi$-divergence estimators (D$\phi$DE). These estimators are also compared with some other methods, including maximum likelihood estimator (MLE), approximate maximum likelihood estimator (AMLE) and estimators based on density power divergence
method (MDPDE).

\noindent Following \cite{Stute1995}, the {K}aplan-{M}eier integral $\int
h(\theta,\alpha)\mathrm{d}\widehat{P}_{n}$  may be written as
\begin{equation*}
    \sum_{i=1}^nW_{in}h(\theta,\alpha,Z_{(i)})
\end{equation*}
where for $1 \leq i \leq n$
\begin{equation*}
  W_{in}=\frac{\delta_{(i)}}{n-i+1}\prod_{j=1}^{i-1}\left[\frac{n-j}{n-j+1}\right]^{\delta_{(j)}}
\end{equation*}
\begin{figure}
\begin{center}
\centerline{\includegraphics[width=8cm]{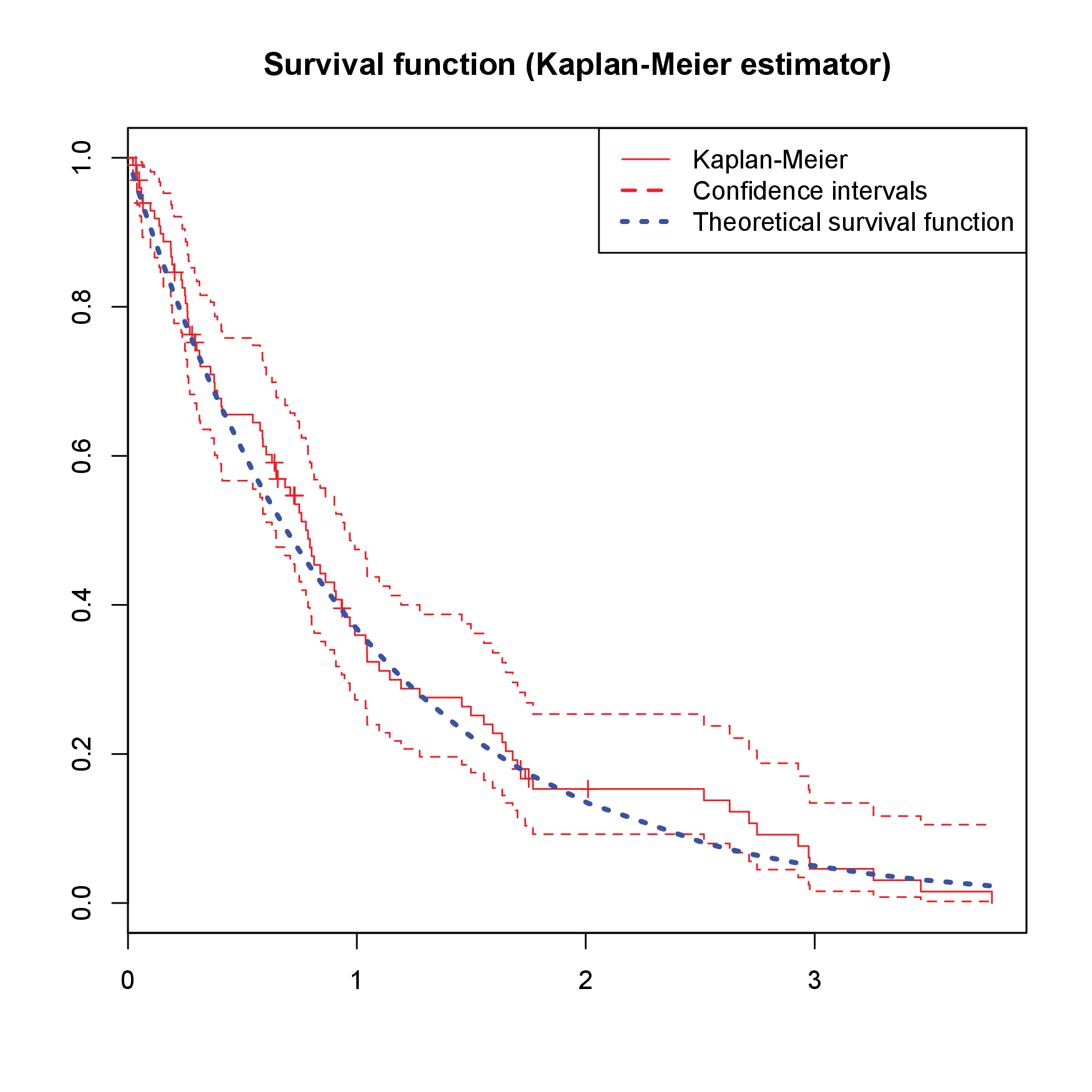}}
\end{center}
\vspace*{-1cm} \hspace*{1.75cm}\begin{minipage} {0.8\linewidth}
\caption{{K}aplan-{M}eier estimator of survival function with confidence intervals.} \label{KMplot}
\end{minipage}
\end{figure}

\noindent Figure \ref{KMplot} presents the {K}aplan-{M}eier estimator of the survival function for a randomly generated exponential sample $\exp(1)$ of size $100$  with $\exp(1/9)$ as censoring distribution. 

\noindent In this simulation study we will use the power divergences family \cite{CressieRead1984}. In this case
\begin{equation*}
\int
h(\theta,\alpha)\mathrm{d}\widehat{P}_{n}=\frac{1}{\gamma-1}\int\left(\frac{p_{\theta}}{p_{\alpha}}\right)^{\gamma-1}~{\rm{d}}P_{\theta}-\frac{1}{\gamma}
\int\left[\left(\frac{p_{\theta}}{p_{\alpha}}\right)^{\gamma}-1\right]~{\rm{d}}\widehat{P}_{n}-\frac{1}{\gamma-1}.
\end{equation*}
\noindent Consider the lifetime distribution to be the one parameter exponential $\exp{(\theta)}$ with density
$p_{\theta}(x)=\theta e^{-\theta x},~x\geq 0$. The MLE of $\theta_0$ is given by
\begin{equation}\label{MLE}
   \widehat{ \theta}_{n,MLE}=\frac{\sum_{i=1}^n\delta_{i}}{\sum_{i=1}^nZ_{i}},
\end{equation}
and the AMLE of \cite{Oakes1986} is defined by
\begin{equation}\label{AMLE}
   \widehat{ \theta}_{n,AMLE}=\frac{\sum_{i=1}^n\delta_{i}}{\sum_{i=1}^nW_{in}Z_{(i)}}.
\end{equation}
It follows that for $\gamma\in\mathbb{R}\setminus\left\{0,1\right\}$
\begin{equation*}
    \frac{1}{\gamma-1}\int\left(\frac{p_{\theta}}{p_{\alpha}}\right)^{\gamma-1}~{\rm{d}}P_{\theta}=\frac{\theta^\gamma\alpha^{1-\gamma}}{(\gamma-1)\left[\gamma\theta+(1-\gamma)\alpha\right]},
\end{equation*}
and
\begin{eqnarray*}
 \int h(\theta,\alpha)\mathrm{d}\widehat{P}_{n}&=&\frac{\theta^\gamma\alpha^{1-\gamma}}{(\gamma-1)\left[\gamma\theta+(1-\gamma)\alpha\right]}\\
&&-\frac{1}{\gamma}\sum_{i=1}^nW_{in}\left[\left(\frac{\theta}{\alpha}\right)^\gamma\exp\left\{-\gamma(\theta-\alpha)Z_{(i)}\right\}-1\right].
\end{eqnarray*}
For $\gamma=0$,
\begin{equation*}
 \int
h(\theta,\alpha)\mathrm{d}\widehat{P}_{n}=\sum_{i=1}^nW_{in}\left[(\theta-\alpha)Z_{(i)}-\log\left(\frac{\theta}{\alpha}\right)\right].
\end{equation*}
Observe that this divergence leads to the AMLE, independently upon the value of $\theta$.

\noindent For $\gamma=1$,
\begin{equation*}
 \int
h(\theta,\alpha)\mathrm{d}\widehat{P}_{n}=\log\left(\frac{\theta}{\alpha}\right)-\frac{(\theta-\alpha)}{\theta}-\sum_{i=1}^nW_{in}\left[ \frac{\theta}{\alpha}\exp\left(-(\theta-\alpha)Z_{(i)}\right)-1
\right].
\end{equation*}

\noindent To make some comparisons, beside dual $\phi$-divergences
estimators, we considered minimum density power divergence
estimators of \cite{BasuBasuJones2006}, (MDPDE's), recall that the density power  divergence between $g$ and another density $f$ is
\begin{equation*}
    d_{\beta}\left(g,f\right)=\int\left\{f^{1+\beta}(z)-\left(1+\frac{1}{\beta}\right)g(z)f^{\beta}(z)+\frac{1}{\beta}g^{1+\beta}(z)\right\}\mathrm{d}z\textrm{ for }\beta>0.
\end{equation*}

\noindent The values of $\gamma$ are chosen to be $-1,~0,~0.5,~1,~2$ which corresponds to the well known standard divergences: $\chi_{m}^{2}-$divergence, $KL_{m}$, the {H}ellinger distance, $KL$  and the $\chi^{2}-$divergence respectively. For the MDPDE's we take the following values of $\beta:~0.1,~0.5,~1$.

\noindent A sample is generated from $\exp(1)$ and $0$, $10$, $25$ of the observations are contaminated by $\exp(5)$ successively. We
have used an exponential censoring scheme, the censoring distribution  is taken to be $\exp(1/9)$, that the proportion of censoring is $10\%$ . The D$\phi$DE's $\widehat{\alpha}_{\phi}(\theta)$ are calculated for samples of sizes $25,~50,~75,~100$ and the hole procedure is repeated $1000$ times. The  value of escort parameter $\theta$ is taken to be the AMLE. We carried out {K}aplan-{M}eier analysis with the \textbf{Survival} package \cite{Survival} within the R Language \cite{Rcore}.

\begin{table}
\caption{MSE of the estimates with $10\%$ of censoring }
\centerline{
\begin{tabular}{lcccccc}
  \hline
  & \multicolumn{6}{c}{$n$} \\
  \cline{2-7}
 & 25 & 50 & 75 & 100 & 150 & 200 \\
  \cline{2-7}
  MLE & 0.0572 & \textbf{0.0250} & \textbf{0.0157} & \textbf{0.0122} & \textbf{0.0079} & \textbf{0.0058} \\
  \cline{2-7}
  $\gamma$&&&&&&\\
  -1 & \textbf{0.0517} & 0.0335 & 0.0188 & 0.0178 & 0.0100 & 0.0090 \\
  0 & 0.0685 & 0.0281 & 0.0166 & 0.0135 & 0.0084 & 0.0062 \\
  0.5 & 0.0727 & 0.0287 & 0.0168 & 0.0138 & 0.0085 & 0.0063 \\
  1 & 0.0824 & 0.0302 & 0.0174 & 0.0143 & 0.0086 & 0.0063 \\
  2 & 0.2533 & 0.1156 & 0.0597 & 0.0436 & 0.0151 & 0.0084 \\
  \cline{2-7}
  $\beta$&&&&&&\\
  0.1 & 0.0643 & 0.0272 & 0.0162 & 0.0131 & 0.0083 & 0.0061 \\
  0.5 & 0.0772 & 0.0368 & 0.0209 & 0.0173 & 0.0112 & 0.0083 \\
  1 & 0.1042 & 0.0506 & 0.0279 & 0.0232 & 0.0154 & 0.0108 \\
   \hline
\end{tabular}
}
\label{MSE1}
\end{table}

\noindent Tables \ref{MSE1} and \ref{MSE2} provide the MSE of various estimates  under the model, according to an an increasing proportion of censoring. As
expected, when there is no contamination, MLE produces most efficient estimators. A close look at the results of the simulations show that the D$\phi$DE's
performs well under the model, when no outliers are generated. For small sample size $n=25$ and $n=50$, the
performance of the estimator under the model is comparable to that of MDPDE's. Indeed in terms of empirical MSE the D$\phi$DE's with $\gamma=-1$ produces a lower MSE than the MDPDE's for all considered values of $\beta$. As $n$ grows up, the MDPDE's prevail.

\begin{table}
\caption{MSE of the estimates with $20\%$ of censoring }
\centerline{
\begin{tabular}{lcccccc}
  \hline
  & \multicolumn{6}{c}{$n$} \\
  \cline{2-7}
 & 25 & 50 & 75 & 100 & 150 & 200 \\
  \cline{2-7}
  MLE & \textbf{0.0627} & \textbf{0.0280} & \textbf{0.0174} & \textbf{0.0134} & \textbf{0.0088} & \textbf{0.0068} \\
  \cline{2-7}
  $\gamma$&&&&&&\\
  -1 & 0.0655 & 0.0395 & 0.0262 & 0.0195 & 0.0154 & 0.0138 \\
  0 & 0.0892 & 0.0395 & 0.0248 & 0.0172 & 0.0113 & 0.0083 \\
  0.5 & 0.0991 & 0.0440 & 0.0273 & 0.0184 & 0.0119 & 0.0087 \\
  1 & 0.1268 & 0.0541 & 0.0336 & 0.0213 & 0.0131 & 0.0094 \\
  2 & 0.3703 & 0.2233 & 0.1919 & 0.1391 & 0.0689 & 0.0510 \\
  \cline{2-7}
  $\beta$&&&&&&\\
  0.1 & 0.0816 & 0.0362 & 0.0224 & 0.0155 & 0.0102 & 0.0075 \\
  0.5 & 0.0919 & 0.0420 & 0.0247 & 0.0171 & 0.0119 & 0.0085 \\
  1 & 0.1166 & 0.0559 & 0.0318 & 0.0218 & 0.0162 & 0.0110 \\
\hline
\end{tabular}
\label{MSE2}
}
\end{table}

\noindent Thus, the D$\phi$DE's are shown to be an attractive alternative to
both the AMLE and MDPDE's in these settings.

\begin{table}
\caption{MSE of the estimates with $20\%$ of  contamination--$10\%$ of censoring}
\centerline{
\begin{tabular}{lcccccc}
  \hline
  & \multicolumn{6}{c}{$n$} \\
  \cline{2-7}
 & 25 & 50 & 75 & 100 & 150 & 200 \\
  \cline{2-7}
  MLE & 0.2413 & 0.1354 & 0.0975 & 0.0916 & 0.0798 & 0.0771 \\
  \cline{2-7}
  $\gamma$&&&&&&\\
     -1 & \textbf{0.0576} & \textbf{0.0617} & \textbf{0.0620} & \textbf{0.0626} & \textbf{0.0605} & \textbf{0.0627} \\
  0 & 0.0852 & 0.0812 & 0.0709 & 0.0710 & 0.0666 & 0.0674 \\
  0.5 & 0.0860 & 0.0820 & 0.0717 & 0.0718 & 0.0676 & 0.0683 \\
  1 & 0.0872 & 0.0826 & 0.0723 & 0.0724 & 0.0682 & 0.0689 \\
  2 & 0.0939 & 0.0843 & 0.0738 & 0.0735 & 0.0692 & 0.0697 \\
  \cline{2-7}
  $\beta$&&&&&&\\
   0.1 & 0.0904 & 0.0905 & 0.0829 & 0.0835 & 0.0834 & 0.0854 \\
  0.5 & 0.1134 & 0.1237 & 0.1243 & 0.1269 & 0.1369 & 0.1405 \\
  1 & 0.1231 & 0.1372 & 0.1424 & 0.1449 & 0.1524 & 0.1547 \\
  \hline
\end{tabular}
\label{MSE3}
}
\end{table}

\begin{table}
\caption{MSE of the estimates with $20\%$ of  contamination--$20\%$ of censoring}
\centerline{
\begin{tabular}{lcccccc}
  \hline
  & \multicolumn{6}{c}{$n$} \\
  \cline{2-7}
 & 25 & 50 & 75 & 100 & 150 & 200 \\
  \cline{2-7}
  MLE & 0.2785 & 0.1629 & 0.1165 & 0.1081 & 0.0962 & 0.0926 \\
  \cline{2-7}
  $\gamma$&&&&&&\\
    -1 & \textbf{0.0624} & \textbf{0.0661} & \textbf{0.0674} & \textbf{0.0684} & \textbf{0.0670} & \textbf{0.0689} \\
  0 & 0.0943 & 0.0898 & 0.0811 & 0.0796 & 0.0751 & 0.0758 \\
  0.5 & 0.0957 & 0.0914 & 0.0826 & 0.0809 & 0.0768 & 0.0774 \\
  1 & 0.0975 & 0.0928 & 0.0840 & 0.0820 & 0.0781 & 0.0784 \\
  2 & 0.1076 & 0.0971 & 0.0872 & 0.0845 & 0.0801 & 0.0801 \\
  \cline{2-7}
  $\beta$&&&&&&\\
  0.1 & 0.0963 & 0.0967 & 0.0891 & 0.0884 & 0.0881 & 0.0900 \\
  0.5 & 0.1127 & 0.1235 & 0.1226 & 0.1241 & 0.1335 & 0.1369 \\
  1 & 0.1225 & 0.1348 & 0.1391 & 0.1409 & 0.1503 & 0.1523 \\
  \hline
\end{tabular}
\label{MSE4}
}
\end{table}

\noindent We now turn to the comparison of these various estimators under contamination. The D$\phi$DE's yield clearly the most robust estimate and outperform the MLE substantially. We can see from Tables \ref{MSE3} and \ref{MSE4}  that the D$\phi$DE with $\gamma=-1$  has the smallest MSE over all other D$\phi$DE's and the MDPDE's for all considered values of $\beta$. As $n$ increases all the D$\phi$DE's compare favorably with MDPE for all  $\beta$.

In the case of long-tailed contamination in the form of an $\exp(0.1)$ distribution, simulations results (not reported in this paper) emphasise that the MDPDE's are more robust than our proposed estimators.

\noindent In conclusion, without contamination the D$\phi$DE's express a good small sample size performance which is comparable to the AMLE and MDPDE's. For medium and large sample sizes the MDPDE's are preferable. Under main body contamination, the D$\phi$DE's are more powerful.

\section{Concluding remarks}
\noindent We have introduced a new estimation procedure in parametric models in the case of right censored data. The method is based on the dual representation of $\phi$-divergences. The estimators are easily
computed and  exhibit appropriate asymptotic behaviour.

We have presented an adaptive choice of the escort parameter $\theta$ that leads to  efficient and robust
estimates. It will be
interesting to investigate theoretically the problem of the choice
of the divergence which leads to an ``\emph{optimal}'' estimate in terms of efficiency and robustness. One approach is to minimize an estimated asymptotic mean squared error of the estimator when it is mathematically tractable, which is not an easy task in the context of censored data and lays beyond the scope of
the present work.

\appendix
\section{Proofs}
\subsection{Proof of Theorem \ref{Theo1}}

\noindent Under the assumptions (R.0), (R.1) and by applying the Strong Law of Large Numbers (SLLN) for censored data, see for instance \cite{StuteWang1993}, \cite{Stute1995} and Proposition 1 in \cite{Wang1999}, we can see that
 \begin{equation}\label{SLLN1}
    \int\frac{\partial}{\partial\alpha}h(\theta,\theta_0)~\mathrm{d}\widehat{P}_{n}\longrightarrow \int\frac{\partial}{\partial\alpha}h(\theta,\theta_0)~\mathrm{d}P_{\theta_0}=0,
 \end{equation}
and
\begin{equation}\label{SLLN2}
    \int\frac{\partial^2}{\partial\alpha\partial\alpha^{\top}}h(\theta,\theta_0)~\mathrm{d}\widehat{P}_{n}\longrightarrow \int\frac{\partial^2}{\partial\alpha\partial\alpha^{\top}}h(\theta,\theta_0)~\mathrm{d}P_{\theta_0}=-S<0,
 \end{equation}

 \noindent Now, for any $\alpha=\theta_0+un^{-1/3}$, with $\|u\|\leq1$, consider a Taylor expansion of $\displaystyle{\int h(\theta ,\alpha )~\mathrm{d}\widehat{P}_{n}}$ in $\alpha$ in a neighborhood of $\theta_0$. Using (R.1), one finds
 \begin{eqnarray}\label{Taylor}
    n\int h(\theta ,\alpha )~\mathrm{d}\widehat{P}_{n}-n\int h(\theta ,\theta_0 )~\mathrm{d}\widehat{P}_{n}&=&n^{2/3}u\int\frac{\partial}{\partial\alpha}h(\theta,\theta_0)~\mathrm{d}\widehat{P}_{n}\\
    \nonumber&&+n^{1/3}\frac{u^2}{2}\int\frac{\partial^2}{\partial\alpha\partial\alpha^{\top}}h(\theta,\theta_0)~\mathrm{d}\widehat{P}_{n}+O(1),
 \end{eqnarray}
uniformly in $u$ with $\|u\|\leq1$. Observe that,
\begin{eqnarray*}
    \left|\int\frac{\partial}{\partial\alpha}h(\theta,\theta_0)~\mathrm{d}\left(\widehat{P}_{n}-P_{\theta_0}\right)\right|&=&\left|\int\left(\widehat{P}_{n}-P_{\theta_0}\right)~\mathrm{d}\left[\frac{\partial}{\partial\alpha}h(\theta,\theta_0)\right]\right|\\
    &\leq&\sup_{x}\left|\widehat{P}_{n}(x)-P_{\theta_0}(x)\right|\int~\mathrm{d}\left|\frac{\partial}{\partial\alpha}h(\theta,\theta_0)\right|.
\end{eqnarray*}
On the other hand, under condition (R.2), by the {LIL} of \cite{FoldesRejto1981}, we have
\begin{equation*}
   \int\frac{\partial}{\partial\alpha}h(\theta,\theta_0)~\mathrm{d}\widehat{P}_{n}=O\left(n^{-1/2}\left(\log\log n\right)^{1/2}\right).
\end{equation*}

\noindent Therefore, using (\ref{SLLN1}) and (\ref{SLLN2}), we obtain for any $\alpha=\theta_0+un^{-1/3}$, with $\|u\|=1$,
\begin{equation*}
    n\int h(\theta ,\alpha )~\mathrm{d}\widehat{P}_{n}-n\int h(\theta ,\theta_0 )~\mathrm{d}\widehat{P}_{n}=O\left(n^{1/6}\left(\log\log n\right)^{1/2}\right)-\frac{1}{2}n^{1/3}S+O(1),
 \end{equation*}

\noindent Observe that the right-hand side vanishes when $\alpha=\theta_0$, and that the left-hand side, by
  (\ref{SLLN2}), becomes negative for all $n$ sufficiently large. Thus, by the continuity of $\alpha\mapsto \int h(\theta ,\alpha )~\mathrm{d}\widehat{P}_{n}$, it holds that as $n\longrightarrow\infty$, with probability one,
  \begin{equation*}
    \alpha\mapsto \int h(\theta ,\alpha )~\mathrm{d}\widehat{P}_{n}
  \end{equation*}
reaches its maximum value at some point $\widehat{\alpha}_\phi(\theta)$ in the interior of  $B(\theta_{0},n^{-1/3})$. Therefore, the estimate $\widehat{\alpha}_\phi(\theta)$ satisfies
 \begin{equation*}
    \int\frac{\partial}{\partial\alpha}h(\theta,\widehat{\alpha}_\phi(\theta))~\mathrm{d}\widehat{P}_{n}=0 \textrm{ and } \|\widehat{\alpha}_\phi(\theta)-\theta_0\|=O(n^{-1/3}).
 \end{equation*}

\subsection{Proof of Theorem \ref{Theo2}}
\noindent Using (R.1), simple calculus give
\begin{equation}
P_{\theta _{0}}\frac{\partial }{\partial \alpha}h(\theta ,\alpha )=0
\label{eqn1}
\end{equation}
and
\begin{equation}
P_{\theta _{0}}\frac{\partial^2}{\partial\alpha\partial\alpha^{\top}}h(\theta ,\theta_0)=-\int \phi ^{\prime \prime }\left(\frac{p_{\theta }}{p_{\theta
_{0}}}\right)\frac{p_{\theta }^{2}}{p_{\theta _{0}}^{3}}\dot{p}_{\theta _{0}}\dot{p}_{\theta _{0}}^{\top}~\mathrm{d}\lambda =:-S.  \label{eqn2}
\end{equation}
Observe that the matrix $S$ is symmetric and positive since the
second derivative $\phi ^{\prime \prime }$ is nonnegative by the
convexity of $\phi $. Let $\displaystyle{U_{n}(\theta _{0}):=\widehat{P}_{n}\frac{\partial }{\partial \alpha}h(\theta ,\theta_0)}$, and use (\ref{eqn1})
and (R.0), (R.3) and (R.4) in connection with the Central Limit Theorem for censored data ({CLT}), see for instance \cite{Stute1995b}, \cite{Wang1999} to
see that
\begin{equation}
\sqrt{n}U_{n}(\theta _{0})\rightarrow \mathcal{N}(0,V).
\label{eqn3}
\end{equation}%
Also, let $\displaystyle{S_{n}(\theta _{0}):=\widehat{P}_{n}\frac{\partial^2}{\partial\alpha\partial\alpha^{\top}}h(\theta ,\theta_0)}$, and use (\ref{eqn2}) and (R.0) in
connection with the {SLLN} to conclude that
\begin{equation}
S_{n}(\theta _{0})\rightarrow -S~~(a.s).  \label{eqn4}
\end{equation}%

\noindent  Using the fact that
$\displaystyle{\widehat{P}_n\frac{\partial }{\partial \alpha}h(\theta,\widehat{\alpha}_{\phi}(\theta))=0}$ and a
Taylor expansion of
$\displaystyle{\widehat{P}_n\frac{\partial }{\partial \alpha}h(\theta,\widehat{\alpha}_{\phi}(\theta))}$ in
$\widehat{\alpha}_{\phi}(\theta)$ around $\theta_0$, we obtain
\begin{equation*}
0=\widehat{P}_n\frac{\partial }{\partial \alpha}h(\theta,\widehat{\alpha}_{\phi}(\theta))=
\widehat{P}_n\frac{\partial }{\partial \alpha}h(\theta,\theta_0)+\left(\widehat{\alpha}_{\phi}(\theta)-\theta_0\right)^{\top}
\widehat{P}_{n}\frac{\partial^2}{\partial\alpha\partial\alpha^{\top}}h(\theta ,\theta_0)+o_P\left(\frac{1}{\sqrt{n}}\right).
\end{equation*}
Hence,
\begin{equation}  \label{eqn6}
\sqrt{n}\left(\widehat{\alpha}_{\phi}(\theta)-\theta_0\right)=
-S_n(\theta_0)^{-1}\sqrt{n} U_n(\theta_0)+o_P(1).
\end{equation}
Using (\ref{eqn3}) and (\ref{eqn4}) and Slutsky Theorem, we
conclude then
\begin{equation}
\sqrt{n}\left(\widehat{\alpha}_{\phi}(\theta)-\theta_0\right)\to \mathcal{N}
\left(0,S^{-1}VS^{-1}\right)
\end{equation}
\subsection{Proof of Theorem \ref{Theo3}}
\noindent By a
Taylor expansion of
$\displaystyle{\widehat{P}_n\frac{\partial }{\partial \alpha}h(\widehat{\theta}_{n},\widehat{\alpha}_{\phi}(\widehat{\theta}_{n}))}$ in
$\widehat{\alpha}_{\phi}(\widehat{\theta}_{n})$ around $\theta_0$, we obtain
\begin{eqnarray*}
0&=&\widehat{P}_n\frac{\partial }{\partial \alpha}h(\widehat{\theta}_{n},\widehat{\alpha}_{\phi}(\theta))=
\widehat{P}_n\frac{\partial }{\partial \alpha}h(\widehat{\theta}_{n},\theta_0)+\left(\widehat{\alpha}_{\phi}(\widehat{\theta}_{n})-\theta_0\right)^{\top}\\
&&\widehat{P}_{n}\frac{\partial^2}{\partial\alpha\partial\alpha^{\top}}h(\widehat{\theta}_{n} ,\theta_0)
+o_P\left(\frac{1}{\sqrt{n}}\right).
\end{eqnarray*}
Taylor expansions of $\displaystyle{\widehat{P}_n\frac{\partial }{\partial \alpha}h(\widehat{\theta}_{n},\widehat{\alpha}_{\phi}(\widehat{\theta}_{n}))}$ and $\displaystyle{\widehat{P}_n\frac{\partial }{\partial \alpha}h(\widehat{\theta}_{n},\theta_0)}$ in
$\widehat{\theta}_{n}$ around $\theta_0$, and the $\sqrt{n}$-consistency of $\widehat{\theta}_{n}$ to $\theta_0$ yield
\begin{eqnarray*}
0&=&\widehat{P}_n\frac{\partial }{\partial \alpha}h(\theta_0,\widehat{\alpha}_{\phi}(\theta))=
\widehat{P}_n\frac{\partial }{\partial \alpha}h(\theta_0,\theta_0)+\left(\widehat{\alpha}_{\phi}(\widehat{\theta}_{n})-\theta_0\right)^{\top}\\
&&\widehat{P}_{n}\frac{\partial^2}{\partial\alpha\partial\alpha^{\top}}h(\widehat{\theta}_{n} ,\theta_0)
+o_P\left(\frac{1}{\sqrt{n}}\right).
\end{eqnarray*}
\noindent Let $\displaystyle{U_{n}:=\widehat{P}_{n}\frac{\partial }{\partial \alpha}h(\theta_0,\theta_0)}$ and $\displaystyle{S_{n}:=\widehat{P}_{n}\frac{\partial }{\partial \alpha}h(\theta_0,\theta_0)}$. By the {CLT}
\begin{equation}
\sqrt{n}U_{n}\rightarrow \mathcal{N}(0,V),
\label{eqn7}
\end{equation}
where $V$ is defined in (\ref{Vnew}).

\noindent Use condition (R.5) and the fact that $\displaystyle{S=-P_{\theta_0}\frac{\partial }{\partial \alpha\partial\alpha^{\top}}h(\theta_0,\theta_0)=-\phi^{\prime\prime}(1)I_{\theta_0}}$,  in connection with Lemma 1 in \cite{Wang1999} to conclude that
\begin{equation}
S_{n}\overset{P}{\longrightarrow} \phi^{\prime\prime}(1)I_{\theta_0}.  \label{eqn8}
\end{equation}%

\noindent The theorem now follows from (\ref{eqn7}), (\ref{eqn8}) and Slutsky's theorem. This concludes the proof.

\end{document}